\documentclass[10pt]{amsart}
\usepackage{fourier}
\usepackage[T1]{fontenc}
\usepackage[utf8]{inputenc}
\usepackage[spanish,english,activeacute]{babel}
\usepackage{mathtools}
\usepackage{enumerate}
\usepackage{hyperref}
\usepackage{stmaryrd}
\usepackage{color} 
\DeclarePairedDelimiter\ceil{\lceil}{\rceil}
\DeclarePairedDelimiter\floor{\lfloor}{\rfloor}

\usepackage{enumerate}

\theoremstyle{plain}
	\newtheorem{thm}{Theorem}[section]
	\newtheorem{prop}[thm]{Proposition}

\theoremstyle{definition}
	\newtheorem{ex}[thm]{Example}
	
\DeclareMathOperator{\Le}{L}

\DeclareMathOperator{\G}{G}
\DeclareMathOperator{\g}{g}
\DeclareMathOperator{\m}{m}

\DeclareMathOperator{\e}{e}

\DeclareMathOperator{\F}{F}
\DeclareMathOperator{\Z}{Z}
\DeclareMathOperator{\PF}{PF}

\DeclareMathOperator{\nc}{nc}
\DeclareMathOperator{\Betti}{Betti}
\DeclareMathOperator{\dis}{d}
\DeclareMathOperator{\cat}{c}

\DeclareMathOperator{\Supp}{Supp}

\title[Numerical semigroups with minimal catenary degree]{Numerical semigroups with embedding dimension three and minimal catenary degree}

\author{Pedro A. García-Sánchez}
\email{pedro@ugr.es}
\address{IEMath-GR and Departamento de Álgebra, Universidad de Granada, 18071 Granada, Spain}
\urladdr{www.ugr.es/local/pedro}

\author{Helena Martín Cruz}
\email{helenamc18@correo.ugr.es}
\address{Departamento de Álgebra, Universidad de Granada, 18071 Granada, Spain}

\begin{document}

\begin{abstract}
We characterize numerical semigroups $S$ with embedding dimension three attaining equality in the inequality $\max\Delta(S)+2\leq \cat(S)$, where $\Delta(S)$ denotes the Delta set of $S$ and $\cat(S)$ denotes the catenary degree of $S$. 
\end{abstract}

\maketitle

\section{Introduction}

It is well known that, for numerical semigroups, the maximum of the Delta set of the semigroup plus two is less than or equal to the catenary degree of the semigroup \cite{GHK}, which in addition is smaller than or equal to the $\omega$-primality \cite{GK}. The $\omega$-primality is smaller than or equal to the tame degree of a numerical semigroup \cite{GHL}. These inequalities relate several important nonunique factorization invariants, and hold for nonfactorial atomic monoids in general. It is also well known that if the presentation of the monoid is generic, then the catenary degree, the $\omega$-primality and the tame degree coincide \cite{STC}. Also, for Krull monoids with finite class group and with prime divisors in all classes, we know that the maximum of the Delta set plus two coincides with the catenary degree in a large number of cases (see \cite[Corollary 4.1]{ggs}), and we still do not know of any Krull monoid where strict inequality holds. In the present paper we show that numerical semigroups show a very different arithmetic behavior.

For numerical semigroups with embedding dimension three, we also know when the catenary degree equals the tame degree \cite{GSV}, or when the tame degree and the $\omega$-pri\-mality coincide \cite{CGSTV}. In this manuscript we focus on the characterization of numerical semigroups with embedding dimension three for which the maximum of the Delta set plus two equals the catenary degree. As an easy exercise, we will also show when this equality holds for numerical semigroups generated by an arithmetic sequence.

It has been shown that the maximum of the Delta set of a numerical semigroup (actually any BF-monoid) is attained in one of its Betti elements \cite{CGSLMS}, and the same holds with the catenary degree \cite{CGSLPR}. Thus the strategy to seek for equality of the maximum of the Delta set plus two and the catenary degree, relies on the study of the sets of factorizations of the Betti elements in our semigroups. For embedding dimension three, these Betti elements are well known (see for instance \cite[Chapter 9]{NS}). For the case of numerical semigroups generated by arithmetic sequences, the result follows easily from the results obtained in \cite{DSNM} and \cite{CT}, for the Delta sets and catenary degrees of these semigroups, respectively.

\section{Preliminaries}\label{pre}

Let $\mathbb{N}$ be the set of nonnegative integers. A numerical semigroup is a nonempty subset of $\mathbb{N}$ that is closed under addition, contains the zero element and whose complement in $\mathbb{N}$ is finite. Every numerical semigroup is finitely generated. If $n_1<\cdots<n_p$ are positive integers with $\gcd(n_1,\dots,n_p)=1$, then the set $\{a_1n_1+\dots+a_pn_p \mid a_1,\dots,a_p \in \mathbb{N}\}$ is a numerical semigroup, say $S$, and every numerical semigroup is of this form. The integers $n_1,\dots,n_p$ are generators of this numerical semigroup, and we write $S=\langle n_1,\dots,n_p\rangle$. If there is no proper subset of $\{n_1,\dots,n_p\}$ that generates $S$, then we say that $S$ is minimally generated by $\{n_1,\dots,n_p\}$. The cardinality of a minimal system of generators of $S$ is called the \emph{embedding dimension} of $S$ and we denote it by $\e(S)$, so with our notation $\e(S)=p$.

From now on, let $S$ be a numerical semigroup minimally generated by $\{n_1,\dots,n_p\}$. The following definitions and results can be found in \cite{NSA} and \cite{NS} in more detail. The homomorphism
\[\pi \colon \mathbb{N}^p \to S \textrm{, } \pi(a_1,\dots,a_p)=\sum_{i=1}^p a_in_i,\]
is the \emph{factorization homomorphism} of $S$, and for $s\in S$, the \emph{set of factorizations} of $s$ in $S$ is the set
$$\mathrm{Z}(s)=\pi^{-1}(s)=\{a \in \mathbb{N}^p \mid \pi(a)=s\}.$$

Let us define the two nonunique factorization invariants on which this work is based. 
The first invariant measures how spread are the lengths of the factorizations of the elements of a numerical semigroup. 

For a factorization $x=(x_1,\dots,x_p)\in \Z(s)$, its \emph{length} is
$$|x|=x_1+\dots+x_p,$$
and the \emph{set of lengths} of factorizations of $s$ is
$$\Le(s)=\{|x| \mid x\in \mathrm{Z}(s)\}.$$
The set of lengths of factorizations of an element in a numerical semigroup is finite. Furthermore, if $S\neq\mathbb{N}$, then there will always be elements with more than one length. Assume that $\mathrm{L}(s)=\{l_1<\cdots<l_k\}$. The \emph{Delta set} of $s$ is the set
$$\Delta(s)=\{l_2-l_1,\dots,l_k-l_{k-1}\}$$
and, if $k=1$, then $\Delta(s)=\emptyset$. The \emph{Delta set} of $S$ is
$$\Delta(S)=\bigcup_{s\in S}\Delta(s).$$

The second invariant we consider in this manuscript, the catenary degree, deals with the distances between factorizations of an element in the numerical semigroup. Let $x=(x_1,\dots,x_p)$, $y=(y_1,\dots,x_p)\in \mathbb{N}^p$ be two factorizations and let
$$x\wedge y=(\min\{x_1,y_1\}, \dots, \min\{x_p,y_p\})$$
be their common part.
The \emph{distance} between $x$ and $y$ is
$$\dis(x,y)=\max\{|x-(x\wedge y)|, |y-(x \wedge y)|\}= \max\{|x|, |y|\}-|x\wedge y|).$$

Let $N\in\mathbb{N}$. A finite sequence of factorizations of $s\in S$, $z_0,\dots,z_k\in \Z(s)$, is called a \emph{$N$-chain of factorizations} if $\dis(z_{i-1},z_i)\leq N$ for all $i\in\{1,\dots,k\}$. The \emph{catenary degree} of $s\in S$, denoted $\cat(s)$, is the least $N\in\mathbb{N}\cup\{\infty\}$ such that for any two factorizations $x$, $y\in \mathrm{Z}(s)$, there is an $N$-chain joining them. The \emph{catenary degree} of $S$, denoted $\cat(S)$, is 
$$\cat(S)=\sup\{\cat(s) \mid s\in S\}.$$

From \cite[Theorem 11]{NSA}, we have the following inequality between the invariants defined above:
$$\max(\Delta(S))+2\leq \cat(S).$$


The complement of $S$ in $\mathbb{N}$ is called the set of \emph{gaps} of $S$, denoted $\G(S)=\mathbb{N}\backslash S$, and its cardinality is called the \emph{genus} of $S$, denoted $\g(S)$. We set $\F(S)=\max(\mathbb{Z}\backslash S)$ and we call it the \emph{Frobenius number} of $S$. The \emph{multiplicity} of $S$ is the smallest nonzero element of $S$ and we denote it by $\m(S)$, so $\m(S)=n_1$.




For a set $X\subseteq \mathbb{N}$, let $X^*$ denote $X\backslash\{0\}$. For a rational number $x$, let $\ceil*{x}$ denote the least integer greater than or equal to $x$, and let $\floor*{x}$ denote the greatest integer less than or equal to $x$.

We say that $x\in \mathbb{Z}$ is a \emph{pseudo-Frobenius number} if $x\notin S$ and $x+s\in S$ for all $s\in S^*$, and denote by $\PF(S)$ the set of these numbers.

A numerical semigroup is \emph{irreducible} if it cannot be expressed as the intersection of two proper over numerical semigroups (with respect to inclusion). It turns out that irreducibility and maximality in the set of numerical semigroups with fixed Frobenius number coincide. We distinguish two cases depending of the parity of the Frobenius number. A numerical semigroup $S$ is \emph{symmetric} (respectively, \emph{pseudo-symmetric}) if $S$ is irreducible and $\F(S)$ is odd (respectively, even).

It is well known (see for instance \cite[Chapter 3]{NS}) that a numerical semigroup $S$ is symmetric if and only if $\g(S)=(\F(S)+1)/2$, and this is equivalent to $\PF(S)=\{\F(S)\}$. Also, $S$ is pseudo-symmetric if and only if $\g(S)=(\F(S)+2)/2$, or equivalently, $\PF(S)=\{\F(S),\F(S)/2\}$. 


A way to construct symmetric numerical semigroups is by gluing symmetric numerical semigroups. Let $S_1$ and $S_2$ be two numerical semigroups minimally generated by $\{n_1,\dots,n_r\}$ and $\{n_{r+1},\dots,n_e\}$, respectively. Let $\lambda\in S_1 \backslash \{n_1,\dots,n_r\}$ and $\mu \in S_2\backslash$ $\{n_{r+1},\dots,n_e\}$ be such that $\gcd(\lambda,\mu)=1$. We say that the numerical semigroup 
$$S=\langle\mu n_1,\dots, \mu n_r, \lambda n_{r+1}, \dots, \lambda n_e\rangle$$
is a \emph{gluing} of $S_1$ and $S_2$. If $S_1$ and $S_2$ are symmetric, then any gluing of $S_1$ and $S_2$ is also symmetric (see for instance \cite[Chapter 8]{NS}).


We define the graph $\nabla_s$ associated to $s\in S$ to be the graph whose set of vertices is $\mathrm{Z}(s)$ and $ab$ is an edge if $a\cdot b\neq 0$ (dot product), that is, there exists $i\in\{1,\dots, p\}$ such that $a_i$, $b_i\neq 0$.

We say that two factorizations $a$ and $b$ of $s$ are $\mathcal R$-related, and we write $a\mathcal{R}b$ or $(a,b)\in\mathcal{R}$, if they belong to the same connected component of $\nabla_s$, that is, there exists a chain of factorizations $a_1,\dots,a_t\in \mathrm{Z}(s)$ such that
\begin{itemize}
\item $a_1=a$, $a_t=b$, and
\item for all $i\in\{1,\dots,t-1\}$, we have $a_i\cdot a_{i+1}\neq 0$.
\end{itemize}

We say that $s\in S$ is a \emph{Betti element} if $\nabla_s$ is not connected. The set of Betti elements of $S$ is denoted by $\Betti(S)$. From these elements one can obtain all minimal presentations of $S$. A presentation for $S$ is a system of generators of the kernel congruence of the factorization homomorphism $\pi$ ($(x,y)\in \ker \pi$ if $\pi(x)=\pi(y)$). A minimal presentation is a presentation that cannot be refined to another presentation, that is, any of its proper subsets is no longer a presentation (see for instance \cite[Chapter 7]{NS}).

In order to obtain a presentation for $S$ we only need, for every $s\in \Betti(S)$ and every connected component $R$, to choose a factorization $x$ and pairs $(x,y)$ such that every two connected components of $\nabla_s$ are connected by a sequence of these factorizations in a way that the pairs of adjacent elements in the sequence are either the ones selected or their symmetry. The least possible number of edges that we need is when we choose the pairs so that we obtain a tree connecting all connected components. Thus, the least possible number of pairs for every $s\in \Betti(S)$ is the number of connected components of $\nabla_s$ minus one \cite[Theorem 5]{NSA}. Thus, the cardinality of any minimal presentation of $S$ equals $\sum_{s\in S} (\nc(\nabla_s)-1)$, where $\nc(\nabla_s)$ is the number of connected components of $\nabla_s$ \cite[Corollary 20]{NSA}, and this cardinality is finite as the set of Betti elements of a numerical semigroup is finite.

Also, the cardinality of a minimal presentation for $S$ is greater than or equal to $\e(S)-1$ \cite[Theorem 8.6]{NS}. It is said that $S$ is a \emph{complete intersection} if the cardinality of any of its minimal presentations equals $\e(S)-1$, that is, that $S$ can be described with the least possible number of relators.

For $x=(x_1,\dots,x_p)\in \mathbb{N}^p$, define $\Supp(x)=\{i\in \{1,\dots,p\} \mid x_i \neq 0\}$. We say that a presentation $\sigma$ of $S$ is \emph{generic} if $\sigma$ is minimal and for all $(x,y)\in \sigma$ we have $\Supp(x+y)=\{1,\dots,p\}$. From \cite[Proposition 5.5]{STC}, we have the uniqueness of generic presentations.





The following results, that are \cite[Theorem 9, Theorem 10]{NSA}, show how minimal presentations, or Betti elements, are key tools to obtain the invariants $\max(\Delta(S))$ and $\cat(S)$.

\begin{thm}\label{thm 16}
Let $S$ be a numerical semigroup. Then
$$\max(\Delta(S))=\max\{\max (\Delta(b)) \mid b \in \Betti(S)\}.$$
\end{thm}

\begin{thm}\label{thm 17}
Let $S$ be a numerical semigroup. Then
$$\cat(S)=\max\{\cat(b) \mid b\in \Betti(S)\}.$$
\end{thm}

Now, we see how the catenary degree can be computed from the factorizations of the Betti elements.

Let $s\in S$ and let $R^s_1,\dots,R^s_{k_s}$ be the different $\mathcal{R}$-classes contained in $\mathrm{Z}(s)$. Set $\mu(s)=\max\big\{r^s_1,\dots,r^s_{k_s}\big\}$, where $r^s_i=\min\left\{|z| \mid z \in R^s_i\right\}$. Define
$$\mu(S)=\max\{\mu(b) \mid b \in \Betti(S)\}.$$

\begin{thm}[{\cite[Theorem 1]{CT}}]\label{thm 20}
Let $S$ be a numerical semigroup. Then
$$\cat(S)=\mu(S).$$
\end{thm}



\begin{ex}\label{ex 18}
Let us see an example of numerical semigroups attaining equality in the inequality $\max(\Delta(S))+2 \leq \cat(S)$. Let $S$ be a numerical semigroup generated by an arithmetic sequence, that is, $S=\langle n, n+k, \dots, n+tk\rangle$ with $1 \leq t < n$ and $\gcd(n,k)=1$. For this type of numerical semigroups, we know the Delta set \cite[Proposition 3.9]{DSNM} and the catenary degree \cite[Theorem 14]{CT}:
$$\Delta(S)=\{k\} \quad \textrm{and} \quad \cat(S)=\ceil*{\frac{n}{t}}+k.$$
For example, if $n=4$, $k=2$ and $t=2$, then $S=\langle 7,9,11,13,15 \rangle$ fulfills $2+2=\max(\Delta(S))+2=\cat(S)=4$.

We have $\max(\Delta(S))+2=\cat(S)$ if and only if  $k+2=\ceil*{\frac{n}{t}}+k$. Equivalently,  $\ceil*{\frac{n}{t}}=2$. If we have three generators ($t=2$), this is equivalent to $n_1 \in \{3,4\}$.
\end{ex}

\section{Numerical semigroups with embedding dimension three}\label{emb dim 3}

We are going to characterize numerical semigroups $S$ with embedding dimension three attaining equality in the inequality $\max(\Delta(S))+2\leq \cat(S)$, to which we will refer as ``$S$ has minimal catenary degree''. The experiments that led to our results were performed using the numerical semigroups package \texttt{numericalsgps} \cite{NSP} for \texttt{GAP} \cite{GAP}. The examples in this section can be reproduced in the following repository\\
\centerline{\url{https://github.com/helenahmc/Examples-min-cat-deg}.}

So let $S=\langle n_1$, $n_2$, $n_3\rangle$, with $2<n_1<n_2<n_3$, be a numerical semigroup with embedding dimension three. Given $\{i, j, k\}=\{1,2,3\}$, define
$$c_i=\min\{k\in\mathbb{Z}^+ \mid kn_i \in \langle n_j,n_k \rangle\}.$$
Then there exist some $r_{ij}$, $r_{ik} \in \mathbb{N}$ such that
\begin{equation}
    c_in_i=r_{ij}n_j+r_{ik}n_k.\label{eq:2}
\end{equation}
From \cite[Example 7.23]{NS}, we know that
\begin{equation}\label{bettis-tres}
    \Betti(S)=\{c_1n_1,c_2n_2,c_3n_3\}.
\end{equation}

Let us give some properties in terms of the parameters already given. First, the\break following theorem is a consequence of some results that Herzog proved in \cite{GR}.

\begin{thm}\label{thm 19}
Let $S$ be a numerical semigroup with embedding dimension three, and let $r_{ij}$ be as defined in \eqref{eq:2}, for all $i$, $j\in\{1,2,3\}$.
\begin{enumerate}
    \item $S$ is a complete intersection if and only if it is symmetric.
    
    \item $S$ is symmetric if and only if there exist $i$, $j \in \{1,2,3\}$ such that $r_{ij}=0$.
    
    \item If $S$ is nonsymmetric, then the integers $r_{ij}$, $r_{ik}$ are positive and unique.
\end{enumerate}
\end{thm}

Note that the uniqueness of the integers $r_{ij}$, $r_{ik}$ in nonsymmetric case can be obtained as a consequence of the uniqueness of generic presentations \cite[Proposition 5.5]{STC}.

From \cite[Lemma 2.3]{DS} and \cite[Lemma 3]{NSEDT}, we have the following result.

\begin{prop}\label{prop 25}
$c_1 > r_{12} + r_{13}$, $c_3 < r_{31} + r_{32}$, and for $i$, $j$, $k\in \{1,2,3\}$,
$$c_i=r_{ji}+r_{ki}.$$
\end{prop}

Let us go back to \eqref{bettis-tres}. Then, $S$ is nonsymmetric if and only if $\#\Betti(S)=3$ (Theorem~\ref{thm 19}) or, equivalently, $S$ has a generic presentation.

In order to compute the maximum of the Delta set and the catenary degree, we distinguish three cases depending on whether $\#\Betti(S)$ is one, two or three. In each case, we will also study when $S$ has minimal catenary degree.

\subsection{A single Betti element}

If $S$ has a single Betti element, that is, $\Betti(S)=\{h=c_1n_1=c_2n_2=c_3n_3\}$, the catenary degree of S is reached in $h$ (Theorem \ref{thm 17}), so $\cat(S)=\cat(h)$. As $S$ is a complete intersection (a minimal presentation has only two relators), $\Z(h)=\{(c_1,0,0),(0,c_2,0),(0,0,c_3)\}$. Thus, $\cat(h)=\max\{c_1,c_2,c_3\}=c_1$ ($n_1<n_2<n_3$ and $c_1n_1=c_2n_2=c_3n_3$ implies $c_1>c_2>c_3$).
\begin{prop}
With the above notation, $\cat(S)=c_1$.
\end{prop}
As $\Le(h)=\{c_3<c_2<c_1\}$, in virtue of Theorem~\ref{thm 16} we have $\max(\Delta(S))=\max(\Delta(h))=\max\{c_2-c_3,c_1-c_2\}$.
\begin{prop}
Under the standing hypothesis, $\max(\Delta(S))=\max\{c_2-c_3,c_1-c_2\}$.
\end{prop}
By \cite[Theorem 12]{ASH}, we have a characterization of this type of numerical semigroups. There exist $p_1>p_2>p_3$ pairwise relatively prime integers greater than one such that $n_1=p_2p_3$, $n_2=p_1p_3$, $n_3=p_1p_2$ and $c_1=p_1$, $c_2=p_2$, $c_3=p_3$.

Let us study when $S$ has minimal catenary degree.

\begin{itemize}
\item Suppose that $\max(\Delta(S))=c_2-c_3$. Then $\max(\Delta(S))+2=\cat(S)$ if and only if $c_1+c_3-c_2=2$. If $c_3=2$, this is a contradiction, because it implies that $c_1=c_2$. Otherwise, it is also a contradiction because the $c_i$ are odd.

\item If $\max(\Delta(S))=c_1-c_2$, then     $\max(\Delta(S))+2=\cat(S)$ if and only if $c_2=2$, which is also a contradiction, because $2=c_2=p_2>p_3=c_3$, and the integers $p_i$ are greater than one.
\end{itemize}

\begin{thm}
Let $S$ be a numerical semigroup with embedding dimension three and a single Betti element. Then
$$\max(\Delta(S))+2<\cat(S).$$
\end{thm}

\subsection{Two Betti elements}

If $S$ has two Betti elements, then $S$ is symmetric. We know exactly how symmetric numerical semigroups with embedding dimension three are thanks to the following result.

\begin{thm}[{\cite[Theorem 9.6]{NS}}]\label{thm 22}
Let $m_1<m_2$ be two relatively prime integers greater than one. Let $a$, $b$ and $c$ be nonnegative integers with $a\geq 2$, $b+c \geq 2$ and $\gcd(a, bm_1 + cm_2)=1$. Then $S=\langle am_1$, $am_2$, $bm_1+cm_2 \rangle$ is a symmetric numerical semigroup with embedding dimension three. Moreover, every symmetric numerical semigroup with embedding dimension three is of this form.
\end{thm}

In virtue of \cite{DSS} and using the notation in Theorem \ref{thm 22}, we have
$$\Betti(S)=\{a(bm_1+cm_2),am_1m_2\},$$
with
\begin{multline*}
   \Z(a(bm_1+cm_2)) =\left\{(b-\floor*{\frac{b}{m_2}}m_2,c+\floor*{\frac{b}{m_2}}m_1,0),\dots,(b-m_2,c+m_1,0),\right.\\
    \left.(b,c,0),(b+m_2,c-m_1,0),\dots,(b+\floor*{\frac{c}{m_1}}m_2,c-\floor*{\frac{c}{m_1}}m_1,0),(0,0,a)\right\},
\end{multline*}
and $\Z(am_1m_2)=\{(m_2,0,0),(0,m_1,0)\}$. Moreover,
\[M_S=\{v= \break(v_1,v_2,v_3)\in \mathbb{Z}^3 \mid v_1n_1+v_2n_2+v_3n_3=0\}\]
is spanned by $\{(m_2,-m_1,0),(b+\lambda m_2,c-\lambda m_1,-a)\}$ for any $\lambda \in \mathbb{Z}$. We choose $\lambda \in\break \left\{-\floor*{\frac{b}{m_2}},\dots,\floor*{\frac{c}{m_1}}\right\}$ such that $|(b+\lambda m_2, c-\lambda m_1, -a)|$ is minimal. We define $$\delta_1=|(m_2,-m_1,0)|=m_2-m_1$$
and
$$\delta_2=\mid|(b+\lambda m_2, c- \lambda m_1, -a)|\mid = \mid b+c+\lambda(m_2-m_1)-a\mid.$$
By \cite[Proposition 5]{DSS}, we have 
\begin{equation}\label{prop 29}
\max(\Delta(S))=\max\{\delta_1,\delta_2\}.
\end{equation}
Using Theorem \ref{thm 20}, we can compute the catenary degree. Indeed, $\cat(S)=\break\max\{\mu(am_1m_2),\mu(a(bm_1+cm_2))\}$. As $m_1<m_2$, $\mu(am_1m_2)=m_2$. The $\mathcal{R}$-classes of $\Z(a(bm_1+cm_2))$ are $\{(0,0,a)\}$ and $\Big\{(b+km_2,c-km_1,0) ~\Big\vert~ k\in\left\{-\floor*{\frac{b}{m_2}},\dots, \floor*{\frac{c}{m_1}}\right\}\Big\}$. Thus, $\mu(a(bm_1+cm_2))=\max\left\{a, b+c-\floor*{\frac{b}{m_2}}(m_2-\right.$ $m_1)\Big\}$. Hence by Theorem~\ref{thm 20},
\begin{equation}\label{prop 30} \cat(S)=\max\left\{m_2,a,b+c-\floor*{\frac{b}{m_2}}(m_2-m_1)\right\}.\end{equation}

We distinguish five cases depending on the position of $a$ in \[\Le(a(bm_1+cm_2))=\{a\}\cup \left\{b+c+k(m_2-m_1) ~\Big\vert~ k\in \left\{-\floor*{\frac{b}{m_2}},\dots,\floor*{\frac{c}{m_1}}\right\}\right\}.\]
\begin{enumerate}[1.]
    \item If $a < b+c-\floor*{\frac{b}{m_2}}(m_2-m_1)- (m_2-m_1)$, then
    \begin{gather*}
    \lambda=-\floor*{\frac{b}{m_2}},  \max(\Delta(S))=b+c-\floor*{\frac{b}{m_2}}(m_2-m_1)-a,\  \textrm{and}\\ \cat(S)=\max\left\{m_2,b+c-\floor*{\frac{b}{m_2}}(m_2-m_1)\right\}.
    \end{gather*}
    \begin{enumerate}[(a)]
        \item If $m_2 < b+c- \floor*{\frac{b}{m_2}}(m_2-m_1)$, then $\cat(S)=b+c-\floor*{\frac{b}{m_2}}(m_2-m_1)$.
        Thus, $\max(\Delta(S))+2=\cat(S)$ if and only if $a=2$.
        
        We have $2=a < b+c-\floor*{\frac{b}{m_2}}(m_2-m_1)- (m_2-m_1)$, that is, $b+c- \floor*{\frac{b}{m_2}}(m_2-m_1) > m_2-m_1+2$, and $b+c-\floor*{\frac{b}{m_2}}(m_2-m_1) > m_2$. As $m_1 \geq 2$, the condition  $b+c-\floor*{\frac{b}{m_2}}(m_2-m_1) > m_2$, implies $a < b+c-\floor*{\frac{b}{m_2}}(m_2-m_1)- (m_2-m_1)$.
        
        To sum up, if $S$ verifies
        $$a=2 \quad \textrm{and} \quad m_2 < b+c-\floor*{\frac{b}{m_2}}(m_2-m_1),$$
        then $S$ has minimal catenary degree.
        For example, if $m_1=5$, $m_2=7$, $b=5$, $c=4$, $a=2$, that is, $S=\langle10,14,53\rangle$, and with the help of \texttt{numericalsgps} we obtain that $\Delta(S)=\{1,2,3,5,7\}$ and $\cat(S)=9$.
        
        \item If $m_2 \geq b+c- \floor*{\frac{b}{m_2}}(m_2-m_1)$, then $\cat(S)=m_2$. Thus $\max(\Delta(S))+2=\cat(S)$ if and only if $m_2=b+c-\floor*{\frac{b}{m_2}}(m_2-m_1)-a+2$. But $m_2 \geq b+c-\floor*{\frac{b}{m_2}}(m_2-m_1)$, so $-a+2 \geq 0$. Hence, $a=2$ and $m_2=b+c-\floor*{\frac{b}{m_2}}(m_2-m_1)$.
        
        We have $2=a < b+c-\floor*{\frac{b}{m_2}}(m_2-m_1)- (m_2-m_1)=m_2-m_2+m_1=m_1$. Thus, $m_1 > 2$.
        
        In this setting, if $S$ verifies
        $$a=2<m_1 \quad \textrm{and} \quad m_2 = b+c-\floor*{\frac{b}{m_2}}(m_2-m_1),$$
        then $S$ has minimal catenary degree.
        For example, if $m_1=3$, $m_2=4$, $b=1$, $c=3$, $a=2$, that is, $S=\langle6,8,15\rangle$, then, $\Delta(S)=\{1,2\}$ and $\cat(S)=4$.
    \end{enumerate}
    
    \item If $b+c-\floor*{\frac{b}{m_2}}(m_2-m_1)- (m_2-m_1) \leq a < b+c-\floor*{\frac{b}{m_2}}(m_2-m_1)$, then
    \begin{gather*}
        \lambda=-\floor*{\frac{b}{m_2}}, \quad \max(\Delta(S))=m_2-m_1 \quad \textrm{and}\\
        \cat(S)=\max\left\{m_2, b+c-\floor*{\frac{b}{m_2}}(m_2-m_1)\right\}.
    \end{gather*}
    
    \begin{enumerate}[(a)]
        \item If $m_2 < b+c-\floor*{\frac{b}{m_2}}(m_2-m_1)$, then $\cat(S)=b+c-\floor*{\frac{b}{m_2}}(m_2-m_1)$.
        Thus, $\max(\Delta(S))+2=\cat(S)$ if and only if $m_2-m_1+2=b+c-\floor*{\frac{b}{m_2}}(m_2-m_1)$, but we have $m_2<b+c-\floor*{\frac{b}{m_2}}(m_2-m_1)$, so $m_1<2$, a contradiction.
        
        \item If $m_2 \geq b+c-\floor*{\frac{b}{m_2}}(m_2-m_1)$, then $\cat(S)=m_2$.
        Thus, $\max(\Delta(S))+2=\cat(S)$ if and only if $m_1=2$.
        
        We have $a < b+c-\floor*{\frac{b}{m_2}}(m_2-2) \leq m_2$. Notice that, as $a\ge 2$, this condition implies that $b+c-\floor*{\frac{b}{m_2}}(m_2-2) -(m_2-2)\le a$.
        
        To sum up, if $S$ verifies 
        $$m_1=2 \quad \textrm{and} \quad  a < b+c-\floor*{\frac{b}{m_2}}(m_2-2) \leq m_2,$$ then $S$ has minimal catenary degree. For example, if $m_1=2$, $m_2=5$, $b=1$, $c=3$, $a=2$, that is, $S=\langle4,10,17\rangle$, then $\Delta(S)=\{1,2,3\}$ and $\cat(S)=5$. 
        
    \end{enumerate}
    
    \item If $b+c-\floor*{\frac{b}{m_2}}(m_2-m_1) \leq a < b+c+\floor*{\frac{c}{m_1}}(m_2-m_1)$, then there exists $k\in\break \left\{-\floor*{\frac{b}{m_2}},\dots,\floor*{\frac{c}{m_1}}-1\right\}$ such that $b+c+k(m_2-m_1) \leq a < b+c+(k+1)(m_2-m_1)$. Hence,
    $$\lambda \textrm{ is either } k \textrm{ or } k+1, \quad \max(\Delta(S))=m_2-m_1 \quad \textrm{and} \quad \cat(S)=\max\{m_2,a\}.$$
    
    \begin{enumerate}[(a)]
        \item If $m_2 < a$, then $\cat(S)=a$.
        
        Thus $\max(\Delta(S))+2=\cat(S)$ if and only if $m_2-m_1+2=a$. As $m_2<a$, $m_1<2$, and we get a contradiction.
        
        \item \label{2.3.b} If $m_2 \geq a$, then $\cat(S)=m_2$.
        
        Thus $\max(\Delta(S))+2=\cat(S)$ if and only if $m_1=2$.
        
        To sum up, if $S$ verifies
        \begin{gather*}
            m_1=2, \quad m_2 \geq a \quad \textrm{and}\\
            b+c-\floor*{\frac{b}{m_2}}(m_2-2) \leq a < b+c+\floor*{\frac{c}{2}}(m_2-2),
        \end{gather*}
        then $S$ has minimal catenary degree. 
        
        
        
As an example, take $m_1=2$, $m_2=7$, $b=1$, $c=2$, $a=5$, and thus, $S=\langle10,16,35\rangle$, and we obtain $\Delta(S)=\{1,2,3,5\}$ and $\cat(S)=7$.
    \end{enumerate}
    
    \item If $b+c+\floor*{\frac{c}{m_1}}(m_2-m_1) \leq a < b+c+\floor*{\frac{c}{m_1}}(m_2-m_1) + m_2-m_1$, then
    $$\lambda=\floor*{\frac{c}{m_1}}, \quad \max(\Delta(S))= m_2-m_1 \quad \textrm{and} \quad \cat(S)=\max\{m_2,a\}.$$
    
    \begin{enumerate}[(a)]
        \item If $m_2 < a$, then $\cat(S)=a$.
        
        Thus $\max(\Delta(S))+2=\cat(S)$ if and only if $m_2-m_1+2=a$. As $m_2<a$, $m_1<2$, and we obtain another contradiction.
        
        \item \label{2.4.b} If $m_2 \geq a$, then $\cat(S)=m_2$.
        
        Thus $\max(\Delta(S))+2=\cat(S)$ if and only if $m_1=2$.
        
        To sum up, if $S$ verifies
        \begin{gather*}
        m_1=2, \quad m_2\geq a \quad \textrm{and}\\
        b+c+\floor*{\frac{c}{2}}(m_2-2) \leq a < b+c+\floor*{\frac{c}{2}}(m_2-2) + m_2-2,
        \end{gather*}
        then $S$ has minimal catenary degree.
        
        
        
        If, for example, $m_1=2$, $m_2=5$, $b=2$, $c=1$, $a=5$, we get  $S=\langle9,10,25\rangle$, and $\Delta(S)=\{1,2,3\}$ and $\cat(S)=5$.
        
    \end{enumerate}
    
    \item If $b+c+\floor*{\frac{c}{m_1}}(m_2-m_1)+m_2-m_1 \leq a$, then
    \begin{gather*}
    \lambda=\floor*{\frac{c}{m_1}}, \quad \max(\Delta(S))= a-b-c-\floor*{\frac{c}{m_1}}(m_2-m_1) \quad \textrm{and}\\
    \cat(S)=\max\{m_2,a\}.
    \end{gather*}
    
    \begin{enumerate}[(a)]
        \item \label{2.5.a} If $m_2 < a$, then $\cat(S)=a$.
        
        Thus $\max(\Delta(S))+2=\cat(S)$ if and only if $b+c+\floor*{\frac{c}{m_1}}(m_2-m_1)=2$.
        
        As $b+c\geq 2$, $m_1<m_2$ and $\floor*{\frac{c}{m_1}} \geq 0$, we deduce that $b+c=2$ and $c<m_1$.
        
        Observe that the conditions $m_2<a$ and $m_1\ge 2$, imply  $2+m_2-m_1=b+c+\floor*{\frac{c}{m_1}}(m_2- m_1) + m_2-m_1 \leq a$-
        
        To sum up, if $S$ verifies
        $$b+c=2, \quad c<m_1 \quad \textrm{and} \quad m_2<a,$$
        then $S$ has minimal catenary degree.
        
        For example, if $m_1=2$, $m_2=3$, $b=1$, $c=1$, $a=7$, that is, $S=\langle5,14,21\rangle$, and we do the computations with \texttt{numericalsgps} we obtain $\Delta(S)=\{1,\ldots,5\}$ and $\cat(S)=7$.
        
        \item \label{2.5.b} If $m_2 \geq a$, then $\cat(S)=m_2$.
        
        Thus $\max(\Delta(S))+2=\cat(S)$ if and only if $a-b-c-\floor*{\frac{c}{m_1}}(m_2-m_1)+2=m_2$.
        
        If $m_2>a$, then $a-b-c-\floor*{\frac{c}{m_1}}(m_2-m_1)+2=m_2$ if and only if $b+c+\floor*{\frac{c}{m_1}}(m_2-m_1) < 2$, in contradiction with $b+c \geq 2$, $m_1 < m_2$ and $\floor*{\frac{c}{m_1}} \geq 0$.
        
        Thus, $m_2=a$, so $b+c+\floor*{\frac{c}{m_1}}(m_2-m_1)=2$ and, as $b+c \geq 2$, $m_1<m_2$ and $\floor*{\frac{c}{m_1}} \geq 0$, $b+c=2$ and $c<m_1$. Then, clearly, $2+m_2-m_1=b+c+\floor*{\frac{c}{m_1}}(m_2-m_1)+m_2-m_1 \leq a$ is equivalent to $m_1 \geq 2$.
        
        To sum up, if $S$ verifies
        $$b+c=2, \quad c<m_1 \quad \textrm{and} \quad m_2=a,$$
        then $S$ has minimal catenary degree. If, for example, $m_1=2$, $m_2=3$, $b=1$, $c=1$, $a=3$, that is, $S=\langle5,6,9\rangle$, we get $\Delta(S)=\{1\}$ and $\cat(S)=3$. 
    \end{enumerate}
\end{enumerate}

Once distinguished the possible cases, note that
\begin{itemize}
    \item we can gather the cases (\ref{2.3.b}) and (\ref{2.4.b}) as $m_1=2$, $m_2\geq a$ and $b+c-\break\floor*{\frac{b}{m_2}}(m_2-2) \leq a < b+c+\floor*{\frac{c}{2}}(m_2-2)+m_2-2$;
    
    \item the same holds for (\ref{2.5.a}) and (\ref{2.5.b}) if we impose $b+c=2$, $c<m_1$ and $m_2\leq a$.
\end{itemize}

Finally, we summarize the results for two Betti elements in the following theorem.

\begin{thm}\label{thm b2}
Let $S$ be a numerical semigroup with embedding dimension three. Assume that $S$ has two Betti elements, and consequently there exist two relatively prime integers $m_1$ and $m_2$ greater than one, with $m_1<m_2$ , and $a$, $b$, and $c$ nonnegative integers with $a\geq 2$, $b+c\geq 2$ and $\gcd(a,bm_1+cm_2)=1$ such that $S=\langle am_1,am_2,bm_1+cm_2 \rangle$. We have
$$\max(\Delta(S))+2=\cat(S)$$
if and only if either
\begin{itemize}
    \item $a=2 \quad \textrm{and} \quad m_2 < b+c-\floor*{\frac{b}{m_2}}(m_2-m_1)$, or
        
    \item $a=2<m_1 \quad \textrm{and} \quad m_2 = b+c-\floor*{\frac{b}{m_2}}(m_2-m_1)$, or
        
    \item $m_1=2 \quad \textrm{and} \quad  a < b+c-\floor*{\frac{b}{m_2}}(m_2-2) \leq m_2$, or
        
    \item $m_1=2, \quad m_2\geq a \quad \textrm{and}$\\ $b+c-\floor*{\frac{b}{m_2}}(m_2-2) \leq a < b+c+\floor*{\frac{c}{2}}(m_2-2)+m_2-2$, or
        
    \item $b+c=2$, $c<m_1$ and $m_2\leq a$.
\end{itemize}
\end{thm}

\begin{ex}
Let us go back to Example \ref{ex 18}, with $n_1=4$ and $t=2$, that is $S=\langle 4,4+k,4+2k\rangle$. We already know that $\max(\Delta(S))+2=\cat(S)=k+2$. Notice $S=\langle4, 4+k, 4+2k\rangle= 2\langle 2,2+k\rangle + (4+k)\mathbb{N}$ is a gluing of $S_1=\langle 2, 2+k\rangle$ and $S_2=\mathbb{N}=\langle 1 \rangle$ by taking $\lambda=4+k \in S_1 \backslash\{2, 2+k\}$ and $\mu=2 \in \mathbb{N}\backslash\{1\}$ ($\gcd(\lambda,\mu)=\gcd(2,k)=1$). Since $S_1$ and $S_2$ are symmetric, so is $S$. A a minimal presentation for $S$ is $\{((2+k,0,0),(0,0,2)),((0,2,0),(1,0,1))\}$ (see \cite[Chapter 8]{NS} to see how a presentation of $S$ is obtained from $S_1$ and $S_2$). Hence $\Betti(S)=\{8+2k,8+4k\}$.
\end{ex}

\subsection{Three Betti elements}

If $S$ has three Betti elements, that is, $S$ is nonsymmetric, then in virtue of \cite{DS}, we set
$$\delta_i=|c_i-r_{ij}-r_{ik}|$$
for every $\{i,j,k\}=\{1,2,3\}$. By Proposition \ref{prop 25}, $\delta_1=c_1-r_{12}-r_{13}$ and $\delta_3=r_{31}+r_{32}-c_3$. In light of \cite{DS}, we have 
\begin{equation}\label{prop 31}
\max(\Delta(S))=\max\{\delta_1,\delta_3\}.
\end{equation}
Also, from \cite[Example 7.23]{NS}, we know that $S$ has a generic presentation. Thus \cite[Corollary 5.8]{STC} is verified, and hence $\cat(S)=\max\{c_1,c_2,c_3,r_{12}+r_{13},r_{21}+r_{23},r_{31}+r_{32}\}$. Moreover, due to Proposition \ref{prop 25}, we obtain the following:
\begin{equation}\label{prop 32}
\cat(S)=\max\{c_1,c_2,r_{21}+r_{23},r_{31}+r_{32}\}.
\end{equation}

We distinguish the cases depending on $\max(\Delta(S))$ equals $\delta_1$ or $\delta_3$, and take into account Proposition \ref{prop 25} (the $r_{ij}$'s are positive).
\begin{enumerate}[1.]
    \item Assume that $c_1-r_{12}-r_{13}>r_{31}+r_{32}-c_3$. This is equivalent to $r_{21}+r_{31}-r_{12}-r_{13}>r_{31}+r_{32}-r_{13}-r_{23}$, which means $r_{21}+r_{23}>c_2$. Hence
    $$\max(\Delta(S))=c_1-r_{12}-r_{13},$$
    and
    $$\cat(S)=\max\{c_1,r_{21}+r_{23},r_{31}+r_{32}\}.$$
    
    \begin{enumerate}[(a)]
        \item \label{3.1.a} If $c_1 \geq \max\{r_{21}+r_{23}, r_{31}+r_{32}\}$, then
        $$\cat(S)=c_1.$$
        
        Thus $\max(\Delta(S))+2=\cat(S)$ if and only if $c_1-r_{12}-r_{13}+2=c_1$. Equivalently, $r_{12}+r_{13}=2$, that is, $r_{12}=r_{13}=1$.
        
        Thus, if $S$ verifies
        $$r_{12}=r_{13}=1, \quad r_{21}+r_{23}>c_2 \quad \textrm{and} \quad c_1 \geq \max\{ r_{21}+r_{23}, r_{31}+r_{32} \},$$
        then $S$ has minimal catenary degree. For example, if $S=\langle4,9,15\rangle$, we have $r_{12}=r_{13}=1$, $3+1=r_{21}+r_{23}>c_2=3$ and $6=c_1 > \max\{ r_{21}+r_{23}=3+1, r_{31}+r_{32}=3+2 \}=5$. With the help of \texttt{numericalsgps} one gets $\Delta(S)=\{1,\ldots, 4\}$ and $\cat(S)=6$.

        \item If $r_{21}+r_{23} \geq \max\{c_1, r_{31}+r_{32}\}$, then
        $$\cat(S)=r_{21}+r_{23}.$$
        
        Thus $\max(\Delta(S))+2=\cat(S)$ if and only if $c_1-r_{12}-r_{13}+2=r_{21}+r_{23}$, that is, $r_{31}-r_{12}-r_{13}-r_{23}+2=0$.
        
        Then we have $r_{21}+r_{23}>c_2$, $r_{21}+r_{23} \geq c_1$, $r_{21}+r_{23} \geq r_{31}+r_{32}$ and $r_{31}-r_{12}-c_3+2=0$.
        
        Observe that $r_{21}+r_{23} \geq c_1=r_{21}+r_{31}$ if and only if $r_{23}\geq r_{31}$, that is, $r_{31}-r_{23}\leq 0$. We then have that $r_{31}-r_{12}-r_{13}-r_{23}+2=0$ is equivalent to $r_{31}=r_{23}$ and $r_{12}=r_{13}=1$. Also,  $r_{21}\ge r_{32}$ is equivalent to $ r_{21}+r_{23}\ge r_{32}+r_{23}=r_{32}+r_{31}$. So these conditions imply $r_{12}+r_{23}\ge \max\{c_1,r_{31}+r_{32}\}$.
        
        Thus, if $S$ verifies
        $$r_{12}=r_{13}=1, \quad r_{21}+r_{23}>c_2, \quad r_{31}= r_{23} \quad \textrm{and} \quad r_{21} \geq r_{32},$$
        then $S$ has minimal catenary degree, and we note that it is a particular case of \ref{3.1.a}.
        
        \item If $r_{31}+r_{32} \geq \max\{c_1, r_{21}+r_{23}\}$, then
        $$\cat(S)=r_{31}+r_{32}.$$
        
        Hence, $\max(\Delta(S))+2=\cat(S)$ if and only if $c_1-r_{12}-r_{13}+2=r_{31}+r_{32}$, that is, $r_{21}-r_{12}-r_{13}-r_{32}+2=0$.
        
        Then we have $r_{21}+r_{23}>c_2$, $r_{31}+r_{32} \geq c_1$, $r_{31}+r_{32} \geq r_{21}+r_{23}$ and $r_{21}-r_{12}-r_{13}-r_{32}+2=0$.
        
        Note that $r_{31}+r_{32} \geq c_1=r_{21}+r_{31}$ is equivalent to $r_{32} \geq r_{21}$, that is, $r_{21}-r_{32} \leq 0$. Then, $r_{21}-r_{12}-r_{13}-r_{32}+2=0$ if and only if  $r_{21}=r_{32}$ and $r_{12}=r_{13}=1$.
        
        Observe that $r_{31}\ge r_{23}$ is equivalent to $ r_{31}+r_{32}\ge r_{23}+r_{32}=r_{23}+r_{21}$. Therefore these conditions imply $r_{31}+r_{32}\ge \max\{c_1,r_{21}+r_{23}\}$.
        
        Thus, if $S$ verifies
        $$r_{12}=r_{13}=1, \quad r_{21} = r_{32}, \quad r_{21}+r_{23} > c_2 \quad \textrm{and} \quad r_{31}\geq r_{23},$$
        then $S$ has minimal catenary degree; again this is a particular instance of \ref{3.1.a}.
    \end{enumerate}
    
    \item If $c_1-r_{12}-r_{13} \leq r_{31}+r_{32}-c_3$, then $r_{21}+r_{31}-r_{12}-r_{13} \leq r_{31}+r_{32}-r_{13}-r_{23}$, equivalently, $r_{21}+r_{23} \leq c_2$. Hence
    $$\max(\Delta(S))=r_{31}+r_{32}-c_3,$$
    and
    $$\cat(S)=\max\{c_1,c_2,r_{31}+r_{32}\}.$$
    
    \begin{enumerate}[(a)]
        \item \label{3.2.a} If $c_1 \geq \max\{c_2, r_{31}+r_{32}\}$, then
        $$\cat(S)=c_1.$$
        
        Thus $\max(\Delta(S))+2=\cat(S)$ if and only if $r_{31}+r_{32}-c_3+2=c_1$, that is, $r_{32}-r_{21}-c_3+2=0$.
        
        Then we have $c_2 \geq r_{21}+r_{23}$, $c_1 \geq c_2$, $c_1 \geq r_{31}+r_{32}$ and $r_{32}-r_{21}-c_3+2=0$.
        
        Observe that $r_{21}+r_{31}=c_1 \geq r_{31}+r_{32}$ if and only if $r_{21} \geq r_{32}$, that is, $r_{32}-r_{21} \leq 0$. The condition $r_{32}-r_{21}-c_3+2=0$ is equivalent to $r_{32}=r_{21}$ and $c_3=2$. Moreover, $c_2 \geq r_{21}+r_{23}=r_{21}+1$ is equivalent to $r_{12}\geq 1$, that is obvious, and $r_{32}=r_{21}$ is equivalent to $r_{32}+r_{31}=r_{21}+r_{31}=c_1$. These conditions imply $c_1 \geq \max\{c_2, r_{31}+r_{32}\}$.
        
        Hence, if $S$ verifies
        $$c_3=2, \quad r_{32}=r_{21} \quad \textrm{and} \quad c_1 \geq c_2,$$
        then $S$ has minimal catenary degree.       If, for example, $S=\langle3,8,13\rangle$, we have $c_3=2$, $1=r_{32}=r_{21}=1$ and $7=c_1>c_2=2$. Also, $\Delta(S)=\{5\}$ and $\cat(S)=7$. 
        
        
        \item \label{3.2.b} If $c_2 \geq \max\{c_1, r_{31}+r_{32}\}$, then
        $$\cat(S)=c_2.$$
        
        Thus $\max(\Delta(S))+2=\cat(S)$ if and only if $r_{31}+r_{32}-c_3+2=c_2$, that is, $r_{31}-r_{12}-c_3+2=0$.
        
        Then we have $c_2 \geq r_{21}+r_{23}$, $c_2 \geq c_1$, $c_2 \geq r_{31}+r_{32}$ and $r_{31}-r_{12}-c_3+2=0$.
        
        Note that $r_{12}+r_{32}=c_2 \geq r_{31}+r_{32}$ is equivalent to $r_{12} \geq r_{31}$, that is, $r_{31}-r_{12} \leq 0$. We then have $r_{31}-r_{12}-c_3+2=0$ if and only if $r_{31}=r_{12}$ and $c_3=2$. Also, the condition $c_2 \geq r_{21}+r_{23}=r_{21}+1$ is equivalent to $c_2 \geq r_{21}+1$. Moreover, note that, if $c_2>c_1$, then the equality $c_2=r_{21}+1$ cannot hold, because this would imply $c_2=r_{21}+1 \leq r_{21}+r_{31}=c_1$, a contradiction.
        
        Observe that $r_{31}+r_{32}=r_{12}+r_{32}=c_2$, so these conditions imply $c_2 \geq \max\{c_1, r_{31}+r_{32}\}$.
        
        Thus, if $S$ verifies
        \begin{gather*}
            c_3=2, \quad r_{31}=r_{12} \quad \textrm{and} \quad c_2 \geq \max \{r_{21}+1, c_1\}\\ \textrm{(excluding the case} \quad c_2>c_1 \quad \textrm{and} \quad c_2 =r_{21}+1 \textrm{),}
        \end{gather*}
        then $S$ has minimal catenary degree.
        For example, if $S=\langle4,5,7\rangle$, we have $c_3=2$, $1=r_{31}=r_{12}=1$, $3=c_2 = \max \{r_{21}+1=2+1, c_1=3\}=3$, $\Delta(S)=\{1\}$ and $\cat(S)=3$.
        
        
        
        \item \label{3.2.c} If $r_{31}+r_{32} \geq \max\{c_1, c_2\}$, then
        $$\cat(S)=r_{31}+r_{32}.$$
        
        Thus $\max(\Delta(S))+2=\cat(S)$ if and only if $r_{31}+r_{32}-c_3+2=r_{31}+r_{32}$, that is, $c_3=2$.
        
        Then we have $c_2 \geq r_{21}+r_{23}=r_{21}+1$, $c_3=2$, $r_{31}+r_{32} \geq c_1$ and $r_{31}+r_{32}\geq c_2$.
        
        As $r_{31}+r_{32}\geq c_1=r_{21}+r_{31}$, then $r_{32}\geq r_{21}$. Thus, $c_2=r_{12}+r_{32} \geq r_{21}+1$ and we can forget about this last condition.
        
        To sum up, if $S$ verifies
        $$c_3=2 \quad \textrm{and} \quad r_{31}+r_{32}\geq \max \{c_1, c_2\},$$
        then $S$ has minimal catenary degree.
        
        
        
        As an example, take $S=\langle7,11,38\rangle$. Then, $c_3=2$, $3+5=r_{31}+r_{32} > c_1=7$, $3+5=r_{31}+r_{32} > c_2=6$, $\Delta(S)=\{1,\ldots, 6\}$ and $\cat(S)=8$.
    \end{enumerate}
    
\end{enumerate}

Once distinguished the possible cases, note that (\ref{3.2.a}) and (\ref{3.2.b}) are particular cases of \ref{3.2.c}. Finally, we gather the results for three Betti elements in the following theorem.

\begin{thm}\label{thm 24}
Let $S$ be a numerical semigroup with embedding dimension three. Assume that $S$ has three Betti elements (equivalently, nonsymmetric). Let $c_i$, $r_{ij}$ be the unique integers fulfilling \eqref{eq:2}. Then 
$$\max(\Delta(S))+2=\cat(S)$$
if and only if either
\begin{itemize}
    \item $r_{12}=r_{13}=1$, $r_{21}+r_{23}>c_2$ and $c_1 \geq \max\{ r_{21}+r_{23}, r_{31}+r_{32}\}$, or
        
    \item $c_3=2$ and $r_{31}+r_{32}\geq \max \{c_1, c_2\}$,
\end{itemize}
\end{thm}

\begin{ex}
Let us revisit once more Example \ref{ex 18}, but now with $n_1=3$ and $t=2$. In this setting $S=\langle 3,3+k,3+2k\rangle$. In light of \cite[Theorem 7]{NSMTF}, $\F(S)=2k$ and $\g(S)=k+1$. As $\g(S)=(\F(S)+2)/2$, $S$ is pseudo-symmetric. Also all generators of $S$ are relatively prime, and thus $S$ has three Betti elements. It is not difficult to see that $c_1=2+k$, and $c_2=c_3=2$.
\end{ex}

\section{Acknowledgements}
The authors would like to thank Alfred Geroldinger for proposing the problem to the second author while she was visiting the University of Graz. 
The first author was partially supported by the Junta de Andaluc\'{\i}a research group FQM-366, and by the project MTM2017-84890-P (MINECO/FEDER, UE). The second author was partially supported by a MINECO collaboration grant in the Departamento de Álgebra, and by the Erasmus+ programme.

\end{document}